# The Lorenz order in graph theory: A new proof and extension of the theorems of Hakimi and of Havel-Hakimi

**Leo Egghe**

Hasselt University, Hasselt, Belgium

E-mail: leo.egghe@uhasselt.be

ORCID: 0000-0001-8419-2932

## Abstract

This paper studies the relation between the Lorenz majorization order and the realizability of degree sequences X of a network in the sense of being graphical or connected graphical (c-graphical) or not. We prove the main result that, if X is dominated (in the Lorenz majorization sense) by X' and X' is (c-) graphical, the X is also (c-) graphical. We present a simple proof and a generalization of the Havel-Hakimi theorem, using the Lorenz order formalism. From this, a classical result of Hakimi on trees follows but also a new generalization to general connected networks. From this, a characterization of c-graphical sequences in terms of the Lorenz majorization order is given.

## 1. Introduction

Let G = (V,E) be a graph or network with V as its vertex set and E as its edge set; #V = N ∈ ℕ, where ℕ is the set of non-null natural numbers {1,2, …}. The degree sequence of this network is denoted as $\Delta_G = (\delta_1, \ldots, \delta_N) \in \mathbb{N}^N$ . All sequences in $\mathbb{N}^N$ are denoted in decreasing order.

Definition: The Lorenz curve (Lorenz, 1905)

Let $X = (x_1, x_2, \ldots, x_N)$ be an N-sequence with $x_j \in \mathbb{R}^+, j = 1, \ldots, N$. If X is an N-sequence, ranked in decreasing order (always used in the sense that ranking is not necessarily strict), then the Lorenz curve of X is the curve in the plane obtained by the line segments connecting the origin (0,0) to the points $\left(\frac{k}{N}, \frac{\sum_{j=1}^{k} x_j}{\sum_{j=1}^{N} x_j}\right)$, k= 1,…,N. For k = N, the endpoint (1,1) is reached.

Definition. The majorization property (Hardy et al., 1934; Marshall et al., 2011).

If X and X′ are N-sequences, ranked in decreasing order, then X is majorized by X′ (equivalently X′ majorizes X), denoted as $X \preccurlyeq_L X'$, if

$$\sum_{j=1}^{k} x_j \leq \sum_{j=1}^{k} x'_j \text{ for } k = 1, \ldots, N-1 \text{ and } \sum_{j=1}^{N} x_j = \sum_{j=1}^{N} x'_j \qquad (1)$$

The index L in $X \preccurlyeq_L X'$ refers to the fact that this order relation corresponds to the order relation between the corresponding Lorenz curves. It is well-known, see e.g., (Marshall et al., 2011, p.14) that $X \preccurlyeq_L Y$ is equivalent to each of the following statements:

(A) $\sum_i \varphi(x_i) \leq \sum_i \varphi(y_i)$ for all continuous, convex functions $\varphi$: $\mathbb{R} \to \mathbb{R}$.

(B) Y can be obtained from X by a finite number of elementary transfers (Muirhead, 1903).

Here an elementary transfer is a transformation from $(x_1, \cdots, x_N)$, where $(x_1, \cdots, x_N)$ is ranked in decreasing order, to $(x_1, \ldots., x_i + h, \ldots., x_j - h, \ldots, x_N)$ where 0 < h ≤ $x_j$.

Definition. Basic transfers (Egghe & Rousseau, 2024)

In the case that the elements in $(x_1, \cdots, x_N)$ are natural numbers, also h can be taken as a natural number, and it can even be taken to be equal to 1. In this case, we will say that this transfer is a basic transfer. It is shown in (Egghe & Rousseau, 2024) how to perform such basic transfers in trees.

We write $X \prec_L Y$ for the strict Lorenz majorization, i.e., $X \preccurlyeq_L Y$ with X ≠ Y.

Definition: Non-normalized Lorenz curves

Let $X = (x_1, x_2, \ldots, x_N)$ be a decreasing N-sequence of non-negative real numbers, then the corresponding non-normalized Lorenz curve is the polygonal line connecting the origin (0,0) with the points $\left(j, \sum_{k=1}^{j} x_j\right)$, j = 1, …, N. This curve ends at the point with coordinates $\left(N, \sum_{k=1}^{N} x_j\right)$.

Definition: The non-normalized (or generalized) majorization order for N-sequences

If X and Y are decreasing N-sequences of non-negative real numbers, then X is majorized by Y, denoted as X ≼ Y if

$$\forall j, j = 1, \dots, N\colon \sum_{k=1}^{j} x_k \le \sum_{k=1}^{j} y_k \tag{2}$$

As for the Lorenz majorization, we write X ≺ Y, for X ≼ Y with X ≠ Y.

Definition (West, 1995). A sequence X = $(x_1, \dots, x_N)$ in $\mathbb{N}^N$ is called graphical if there exists a network G such that $\Delta_G = X$. Such a network is said to realize X.

If the sequence X can be realized through a connected network we say that X is c-graphical, and X is c-realized by this network.

The next theorem provides a sufficient condition for X to be graphical or c-graphical.

Theorem 1

Let X = $(x_j)_{j=1,\dots,N}$, X′ = $(x_j')_{j=1,\dots,N}$ ∈ $\mathbb{N}^N$, with X ≼ X′ = $\Delta_G$ (hence X′ is graphical) and $\sum_j x_j = \sum_j x_j' = \sum_j \delta_j$. Then also X is graphical. Moreover, if X′ is c-graphical then also X is c-graphical.

Proof. Assume that X ≼ X′ = $\Delta_G$ (G connected) and $\sum_j x_j = \sum_j x_j' = \sum_j \delta_j$ . As $\sum_j x_j = \sum_j x_j'$ the generalized majorization order coincides with the standard majorization order, and hence we may apply the discrete version of Muirhead's theorem. Hence we know that we can obtain X′ from a finite number of basic transfers. From now on we continue by induction. Assume first that X′ can be reached from X through one basic transfer. Then, with i < j, X′ = $\Delta_G$ = $(x_1, \dots, x_{i-1}, x_i + 1, x_{i+1}, \dots, x_{j-1}, x_j - 1, x_{j+1}, \dots, x_N) = (\delta'_1, \dots, \delta'_N)$, where $\Delta_G$ is the degree sequence of the connected network G. Now X = $(x_1, \dots, x_N)$ can be rewritten as: $(\delta'_1, \dots, \delta'_{i-1}, \delta'_i - 1, \delta'_{i+1}, \dots, \delta'_{j-1}, \delta'_j + 1, \delta'_{j+1}, \dots, \delta'_N)$. The network we are looking for, which realizes X has, of course, the same nodes as G. X is then c-graphical if there exists a point k, different from the points at ranks i and j, in the network G that has a direct link with i and not with j. Then we delete the link between i and k and add a link between j and k, making sure that the network stays connected. Such a point k ≠ i,j exists. Indeed, as the network G is connected there exists a path $P_{ij}$ connecting i to j which stays invariant. Then, in addition to the path $P_{ij}$, we have for i and j a number of degrees equal to $\delta'_i - 1$, respectively $\delta'_j + 1$ in the given network G. As $x_i \geq x_j$ we have $\delta'_i - 1 \geq \delta'_j + 1 > \delta'_j > \delta'_j - 1$ and hence there exists k ≠ i,j such that k and i are linked, k and j are not linked, and k ∉ $P_{ij}$. Hence if we remove a link between k and i and add a link between j and k we obtain a connected network with X = Δ, which is connected because the path $P_{ij} \cup \{(j.k)\}$ still exists and no other links have been deleted.

Assume now that n basic transformations are needed to go from X to X′ and that we already know that the theorem holds for the case of n-1 basic transformations. If X′ is obtained from X in n basic transformations, then we know, by the induction base, that also the sequence obtained after n-1 basic transformations is c-graphical. Then, it follows from the induction hypothesis that also X is c-graphical.

If X′= $\Delta_G$ is graphical (with G not connected) then it follows from the above that X = Δ is graphical. Then any point k linked to i and not linked to point j can be used (and such a k exists as $\delta_i' > \delta_j'$). □

The opposites of theorem 1 do not hold: if X is c-graphical then X′ is not necessarily c-graphical. Further, the theorem does not hold if $\sum_j x_j \neq \sum_j x_j'$, even if both sums are even.

To show that X ≼ X′, $\sum_j x_j = \sum_j x_j'$ and X c-graphical, does not imply that X′ is graphical, it suffices to provide one example. Consider X = (4,3,3,3,1). This sequence (with N = 5) is c-graphical as it is the degree sequence of the network shown in Fig. 1.

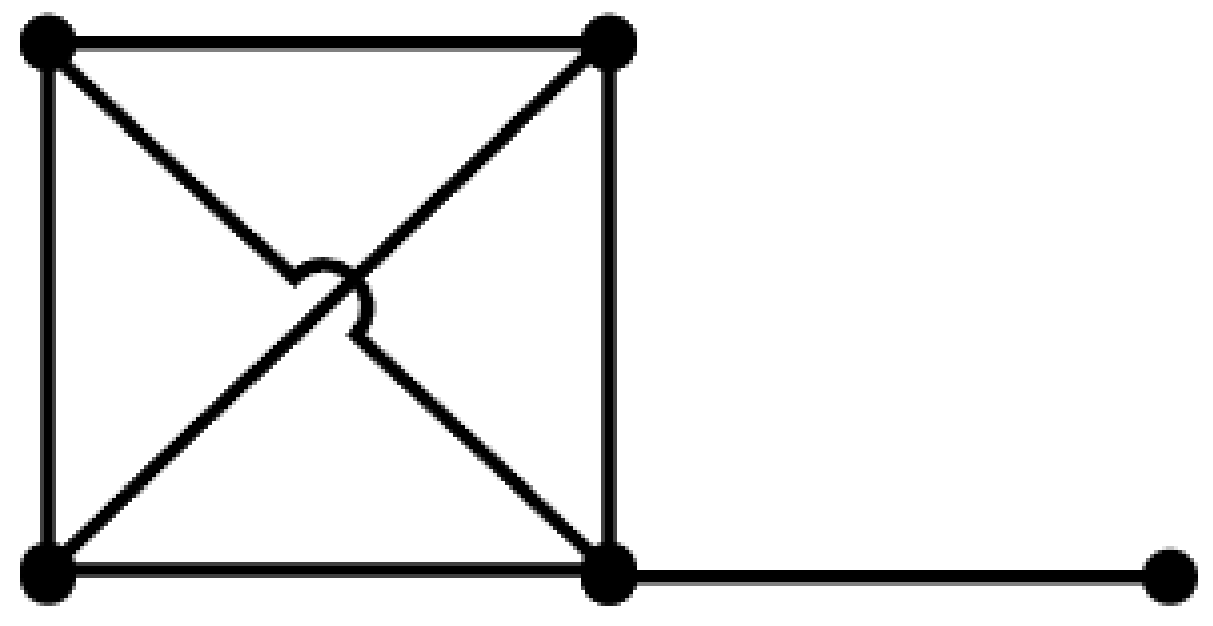

Fig. 1

Now, X = (4,3,3,3,1) ≼ (5,3,3,2,1) = X′. The sequences X and X′ have an equal sum, but X′ is not graphical as there cannot be a node with degree 5 in a 5-node simple network. Even if all nodes have a degree at most 4, the theorem does not hold. Indeed, consider X = (4,4,4,1,1). This sequence is not graphical (three nodes with degree 4 imply that all other nodes have at least degree 3) and (4,3,3,3,1) ≼ (4,4,4,1,1).

Theorem 1 is false if $\sum_j x_j \neq \sum_j x_j'$. Indeed let X′ = (4,4,4,4,4) which is c-graphical as it is the degree sequence of the complete 5-network. Now, X = (4,4,2,1,1) is not graphical (two nodes with degree 4 imply that all other nodes have at least degree 2). Yet, (4,4,2,1,1) ≼ (4,4,4,4,4), proving that the theorem does not hold.

Theorem 1 also does not hold for the classical Lorenz majorization.

Proposition 1

(i) If $X \preccurlyeq_L X'$, X is c-graphical and $\sum_j x_j'$ is even, then X′ is not necessarily graphical.

(ii) If $X \preccurlyeq_L X'$, X′ is c-graphical and $\sum_j x_j$ is even, then X is not necessarily c-graphical.

Proof. It suffices to give one counterexample for each case.

Case (i). Consider X = (4,4,4,4,4) which is obviously c-graphical and X′ = (4,4,2,1,1). We see that $X \preccurlyeq_L X'$,$\sum_j x_j'$ is even and we already know that X′ is not graphical.

Case (ii). Consider X = (2,1,1,1,1) $\preccurlyeq_L X'$ = (4,1,1,1,1). X′ is the delta sequence of a 5-star and hence is c-graphical; $\sum_j x_j$ is even, yet X = (2,1,1,1,1) is not c-graphical (but it is graphical).

Next we present a new proof of the graphical case of Theorem 1.

We recall the famous Erdös-Gallai (1960) theorem which states that

X = $(x_1, \ldots, x_N) \in \mathbb{N}^N$, ordered decreasingly, is graphical if and only if $\sum_j x_j$ is even and $\forall k = 1, \ldots, N$: $\sum_{i=1}^{k} x_i \le k(k-1) + \sum_{i=k+1}^{N} \min(x_i, k)$, where $\sum_{i=N+1}^{N} \min(x_i, k) = 0$.

Theorem 1A. Let X = $(x_j)_{j=1,\ldots,N}$, X' = $(x'_j)_{j=1,\ldots,N} \in \mathbb{N}^N$, with X $\preccurlyeq$ X' and $\sum_j x_j = \sum_j x'_j$, then X' graphical implies that also X is graphical.

We first prove the following Lemma.

Lemma. If X = $(x_j)_{j=1,\ldots,N} \preccurlyeq$ X' = $(x'_j)_{j=1,\ldots,N}$ and $\sum_j x_j = \sum_j x'_j$, then

$$\forall k = 1, \ldots, N-1: \sum_{j=k+1}^{N} \min(x_j, k) \ \ge\ \sum_{j=k+1}^{N} \min(x'_j, k) \qquad (3)$$

Proof. We apply Muirhead's theorem and if X' can be obtained from X through one basic transformation there exist i and j such that

$$X = (x_1, \ldots, x_N) \preccurlyeq X' = (x'_1, \ldots, x'_N) =$$
$$(x_1, \ldots, x_{i-1}, x_i + 1, x_{i+1}, \ldots, x_{j-1}, x_j - 1, x_{j+1}, \ldots, x_N) \qquad (4)$$

Hence: $\forall p = 1, \ldots, N,\ p \ne i, j:\ x_p = x'_p$. We also see that $\forall k = 1, \ldots, N-1$: $\sum_{l=k+1}^{N} x_l \ \ge\ \sum_{l=k+1}^{N} x'_l$.

It follows that

$$\forall\, k = j, \ldots, N-1: \sum_{p=k+1}^{N} \min(x_p, k) = \sum_{p=k+1}^{N} \min(x'_p, k)$$

$$\forall\, k = i, \dots, j-1{:}\ \sum_{p=k+1}^{N} \min(x_p, k) = \min(x_j, k) + \sum_{\substack{p=k+1 \\ p\neq j}}^{N} \min(x_p, k)$$

$$\geq \min(x_{j-1}, k)$$

$$+ \sum_{\substack{p=k+1 \\ p\neq j}}^{N} \min(x'_p, k)$$

$$= \min(x'_j, k) + \sum_{\substack{p=k+1 \\ p\neq j}}^{N} \min(x'_p, k) = \sum_{p=k+1}^{N} min(x'_p, k)$$

$$\forall\, k = 1, \dots, i-1{:}\ \sum_{p=k+1}^{N} \min(x_p, k)$$

$$= \min(x_i, k) + \min(x_j, k) + \sum_{\substack{p=k+1 \\ p\neq i,j}}^{N} \min(x_p, k)$$

$$= \min(x_i, k) + \min(x_j, k) + \sum_{\substack{p=k+1 \\ p\neq i,j}}^{N} \min(x'_p, k)\ .$$

Now, $min(x_i, k) + min(x_j, k) \geq (*) min(x_i + 1, k) + min(x_j - 1, k) = min(x'_i, k) + min(x'_j, k)$. The inequality (*) follows from the facts that if $k \leq x_{j-1}$ then the inequality (*) becomes 2k≥2k, which is correct; if $x_j \leq k \leq x_i$ then the inequality (*) becomes $k + x_j \geq k + x_j - 1$ which is correct; and finally, if $k \geq x_i + i$ then (*) becomes $x_i + x_j \geq x_i + 1 + x_j - 1$ which is also correct.

In this way, the lemma, for one basic transformation, is shown in each case. It follows by induction that the lemma holds for an arbitrary number of basic transformations. □

Proof of Theorem 1A.

$X \preccurlyeq X'$ implies that $\forall k = 1, \dots, N$: $\sum_{j=1}^{k} x_j \le \sum_{j=1}^{k} x'_j$ . If now X′ is graphical, it follows from the Erdös-Gallai theorem that, for k = 1, …, N:

$$\sum_{j=1}^{k} x'_j \le k(k-1) + \sum_{j=k+1}^{N} min(x'_j, k) \qquad (5)$$

Continuing, using the lemma, we obtain that (5) $\le k(k-1) + \sum_{j=k+1}^{N} \min(x_j, k)$ . Then, applying Erdös-Gallai again, we find that X is graphical.□

We recall the Havel-Hakimi theorem.

Theorem (Havel-Hakimi)

X = $(x_1, \dots, x_N) \in \mathbb{N}^N$ is graphical if and only if $X_H = (x_2 - 1, \dots, x_{x_1+1} - 1, x_{x_1+2}, x_{x_1+3}, \dots, x_N) \in \mathbb{N}^{N-1}$ is graphical.

Note that the Havel-Hakimi theorem does not say anything about being c-graphical.

As a corollary of the previous result, we give a new proof of the Havel-Hakimi theorem.

Proof. Since the other direction is trivial we only present the part that if $X = (x_1, \dots, x_N)$ is graphical, then also $X_H = (x_2 - 1, \dots, x_{x_1+1} - 1, x_{x_1+2}, x_{x_1+3}, \dots, x_N) \in \mathbb{N}^{N-1}$ is graphical.

We take X in decreasing order. We remove from the graph associated with X the first vertex and the corresponding $x_1$ links. The resulting sequence X′ is then also graphical. Hence, there exist vertices $i_1, \dots, i_{x_1}$ with $2 \le i_1 < \cdots < i_{x_1} \le N$ such that $X' = (x_2 - \varepsilon_2, \dots, x_N - \varepsilon_N)$ where $\varepsilon_j = 1 \text{ for } j = i_1, \dots, i_{x_1}$ and $\varepsilon_j = 0 \text{ for } j \neq$ $i_1, \ldots, i_{x_1}$. Since $2 \le i_1, 3 \le i_2, \ldots, x_1 + 1 \le i_{x_1}$ (and X is decreasing) we have that $X_H \preccurlyeq X$ with equal sums (since X results from $X_H$ by a finite number ($\le x_1$) of elementary transfers from $i \ge x_1 + 2$ to $i \le x_1 + 1$). We conclude by Theorem 1 that $X_H$ is graphical. □

The following generalization of Havel-Hakimi is due to (Wang & Kleitman, 1973) and can also be found in (West, 2001, p.52 as Exercise 1.3.58)

Theorem (Wang and Kleitman, 1973)

$\forall\, k \in \{1, \ldots, N\} : X = (x_1, \ldots, x_N) = \Delta_G$ is graphical if and only if $X'_k$ is graphical, where $X'_k$ is the degree sequence of the graph obtained from G by deleting node k and subtracting 1 from the $x_k$ largest remaining values in X.

Proof. We only show the non-trivial direction. So, we have $X = (x_1, \ldots, x_N)$, ranked decreasingly and being graphical, and we have to show that $X'_k = (x_1 - \gamma_1, \ldots, x_{k-1} - \gamma_k, x_{k+1} - \gamma_{k+1}, \ldots, x_N - \gamma_N)$ is graphical, where $\gamma_i = 1$ if i belongs to the first $x_k$ indices, i≠k, and $\gamma_i = 0$ otherwise.

Indeed, since *X* is graphical, also $X_k$ obtained from *X* by deleting vertex *k* and its $x_k$ links, is graphical. Hence there exist indices $i_1, \ldots, i_{x_k}$ with $1 \le i_1 < \cdots < i_{x_k} \le N$ , such that $X'_k = \left(x_j - \gamma_j\right)_{j=1}^{N}$ (j ≠k) with $\gamma_i = 1$ if i belongs to the first $x_k$ indices, i≠k, and $\gamma_i = 0$ otherwise. Hence $X'_k \preccurlyeq X_k$, where both sequences have the same sum. Hence by Theorem 1, $X'_k$ is graphical. □

The theorem of Wang and Kleitman, and hence also the theorem of Havel-Hakimi, can further be generalized as follows.

Theorem 2. Let $X = (x_1, \dots, x_N)$ be ranked decreasingly. Then we have for every $k \in \{1, \dots, N\}$ and every $n \in \{1, \dots, x_k\}$ that

$$X \text{ is graphical} \Leftrightarrow X'_{n,k} \text{ is graphical}$$

where $X'_{n,k}$ is obtained from $X$ by changing $x_k$ to $x_k - n$ and subtracting 1 from the n largest remaining values in $X$.

Proof. For the sake of completeness, we present the complete proof, which is similar to the preceding ones

Part 1: the (trivial) arrow from right to left. If, for n and k as defined above, $X'_{n,k}$ is graphical, then adding to te corresponding graph n links from vertex k to the first n vertices (excluding k) yields a graph for which $X$ is the degree sequence. This shows that $X$ is graphical.

Part 2: the arrow from left to right. Suppose now that $X$ is graphical. Hence, for n and k as defined above, we can delete n links from  vertex k, yielding a graph with degree sequence $X_{n,k}$ . Now, since X is decreasing and by the definition of $X'_{n,k}$ we have $X'_{n,k} \preccurlyeq X_{n,k}$ . By Theorem 1 this shows that $X'_{n,k}$ is graphical. □

Remark. If $n = x_k$ it is irrelevant whether we keep vertex k as an isolated point or delete the vertex k (as in the Havel-Hakimi and Wang-Kleitman results).

An example. Consider the 4-chain (N=5) with Δ = (2,2,2,1,1). Let k = 3, $x_k$ = 2 , n=1, then $X'_{n,k}$ =(2,1,1,1,1) as illustrated in Fig. 2.

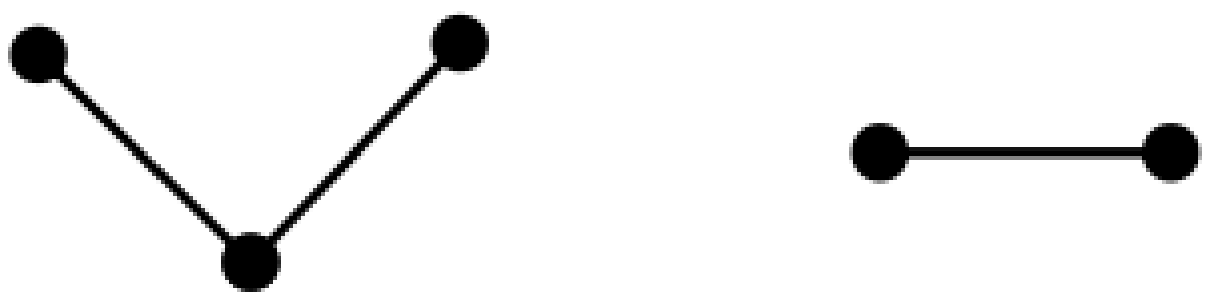

Fig. 2

Application

My new criterion for realizability often gives a shorter way for investigating whether a sequence is graphical or not. When applying Havel-Hakimi or Wang-Kleitman one must apply their criterion until a zero sequence is reached. The new criterion often can be applied until a constant sequence $A = \left(\underbrace{a, \dots, a}_{N\ times}\right)$, with *a* a natural number smaller than or equal to N-1, and *Na* even, is reached. Indeed, A ≼ Δ(G), with Δ(G) the degree sequence of every graph G with sum *Na*, see also West (2001, exercise 1.3.27).

Illustration

To investigate if (5,4,4,3,3,3) is graphical we check with Havel-Hakimi: (5,4,4,3,3,3) → (3,3,2,2,2) → (2,2,1,1) → (1,1,0) → (0,0).

With our new criterion, this becomes (with k=1, n=2): (5,4,4,3,3,3) → (3,3,3,3,3,3), which is finished in fewer steps.

In both cases we conclude that (5,4,4,3,3,3) is graphical, see Fig. 3.

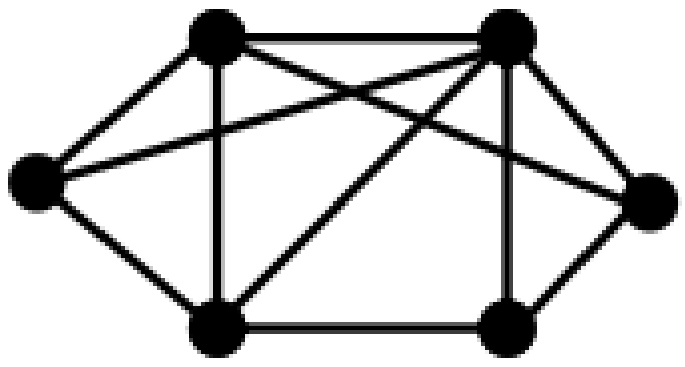

Fig.3. Network with degree sequence (5,4,4,3,3,3)

Next, we show that Hakimi's theorem follows from Theorem 2.

Theorem (Hakimi, 1962). Every sequence $X \in \mathbb{N}^N$ with $\sum_{j=1}^{N} x_j = 2(N-1)$ is the degree sequence of a tree.

Proof. As $X \in \mathbb{N}^N$, every $x_j \geq 1$. Then the degree sequence of a star S $\Delta_S = \left(N-1, \underbrace{1,\dots,1}_{N-1\ times}\right)$ is the largest sequence for the generalized majorization order among all sequences with $\sum_{j=1}^{N} x_j = 2(N-1)$: $X \preccurlyeq \Delta_S$. It follows from Theorem 1 that X is c-graphical (as the star is). Now, as the network realizing X is connected and $\sum_{j=1}^{N} x_j = 2(N-1)$, this network is a tree.

Next, we will generalize Hakimi's theorem to general connected networks. Hakimi's theorem deals with trees, i.e. connected networks such that the sum of its degrees is $\sum_{j=1}^{N} x_j = 2(N-1)$. For a general connected network, we have: 2(N-1) $\leq \sum_{j=1}^{N} x_j \leq N(N-1)$, with $\sum_{j=1}^{N} x_j$ even. Hence the degree sum can be written as 2(N-1) + 2d, with d = 0, …, (N-1)(N-2)/2.

Notation. For a given natural number N > 0, and each d in the set $\left\{0,1,\dots,\frac{(N-1)(N-2)}{2}\right\}$ we denote by $T^d(N)$, or simply $T^d$ if N is known, the set of all connected networks with $\sum_{j=1}^{N} x_j = 2(N-1) + 2d$ .

We recall that a partition of a set is an arrangement of its elements into non-empty subsets, in such a way that every element is included in exactly one subset.

Lemma.

$\forall N \in \mathbb{N}$, $\{T^d(N)\}$ is a partition of the set of all connected networks with N nodes.

Proof. As 2(N-1) + (N-1)(N-2) = N(N-1). It suffices to show that for all d in $\left\{0,1,\dots,\frac{(N-1)(N-2)}{2}\right\}$, $T^d(N) \neq \emptyset$. This is clear as $T^0(N)$ consists of all trees with N nodes. Then $T^d(N)$ consists of all trees in which d extra pairs of nodes are connected.

Construction

Let the natural number N be given and denote the N-star by $S_0$. Then its degree sequence is

$$\Delta_{S_0} = \left(N-1,\ \underbrace{1,\dots,1}_{(N-1)times}\right) \tag{6}$$

Next, we start the first cycle in this construction by choosing one of the nodes with degree 1 and connecting it, step by step (hence in N-2 steps) to each of the other nodes with degree one. The connected graph obtained after d steps, i.e. with d edges added, is denoted as $S_d$. Then

$$\Delta_{S_1} = \left( N-1, 2{,}2, \underbrace{1, \dots, 1}_{(N-3)times} \right) \tag{7}$$

$$\Delta_{S_2} = \left( N-1, 3{,}2{,}2, \underbrace{1, \dots, 1}_{(N-4)times} \right) \tag{8}$$

$$\Delta_{S_3} = \left( N-1, 4{,}2{,}2{,}2 \quad \underbrace{1, \dots, 1}_{(N-3)times} \right) \tag{9}$$

ending the first cycle with

$$\Delta_{S_{N-2}} = \left( N-1, N-1, \underbrace{2, \dots, 2}_{(N-2)times} \right) \tag{10}$$

Now, the second cycle starts by choosing one of the nodes with degree 2 and connecting it in N-3 steps to each of the other nodes with degree 2. This gives:

$$\Delta_{S_{N-1}} = \left( N-1, N-1{,}3{,}3, \underbrace{2, \dots, 2}_{(N-4)times} \right) \tag{11}$$

$$\Delta_{S_N} = \left( N-1, N-1{,}4{,}3{,}3 \quad \underbrace{2, \dots, 2}_{(N-5)times} \right) \tag{12}$$

ending with

$$\Delta_{S_{2N-5}} = \left( N-1, N-1, N-1, \underbrace{3, \ldots, 3}_{(N-3) times} \right) \qquad (13)$$

This procedure continues until we reach

$$\Delta_{S_{N(N-3)/2}} = \left( \underbrace{N-1, \ldots, N-1}_{(N-2) times}, N-2, N-2 \right) \qquad (14)$$

ending the construction with the complete N-graph:

$$\Delta_{S_{(N-2)(N-1)/2}} = \left( \underbrace{N-1, \ldots, N-1}_{N\ times} \right) \qquad (15)$$

Recall that we are allowed to choose nodes freely as isomorphic networks are considered to be the same. This construction has N-1 cycles. We observe that, by construction, for fixed N and fixed d, d = 0, …, (N-1)(N-2)/2, $\Delta_{S_d}$ is the highest in the Lorenz order $\preccurlyeq$ .

Notation. We recall that in the construction above, d denotes the number of links added to the N-star. Hence, $d \in \left\{0, 1, \ldots, \frac{(N-1)(N-2)}{2}\right\}$. For each d, we denote by i(d) the number of nodes with degree (N-1). Hence i(0) = … = i(N-3) = 1; i(N-2) = 2, and so on, ending with $i\left(\frac{(N-1)(N-2)}{2}\right) = N$.

We see that each cycle ends when one more node has a degree (N-1). This means that at the end of the cycle i(d) increases by 1. At that point, the other coordinates are equal to the new i(d) value. We further see that the number of links to be added to obtain a new node with degree (N-1) decreases in each cycle. First, it is (N-2), then (N-3), ending with just one link to be added.

We also have the following formula where j(d) is equal to the number of steps taken after the latest increase of i(d). It is also equal to the value of component (i(d) + 1) minus the value i(d):

$$d = \sum_{k=2}^{i(d)}(N-k) + j(d) \tag{16}$$

with $i(d) \in \{1, \dots, N-1\}$, $j(d) \in \{0, \dots, N-i(d)-2\}$, $\sum_{k=2}^{1}(N-k) = 0$ and $j(N-1) = 0$. From the construction and definitions above we see that

$$\Delta_{S_d} = \Bigg(\underbrace{N-1, \dots, N-1}_{i(d) times}, i(d) + j(d), \underbrace{i(d)+1, \dots, i(d)+1}_{j(d) times}, \underbrace{i(d), \dots, i(d)}_{(N-i(d)-j(d)-1) times}\Bigg) \tag{17}$$

Using this construction we can formulate the following extension of Hakimi (1962).

Proposition 2

Let N > 0 be given, let $X \in \mathbb{N}^N$ and let $\sum_{j=1}^{N} x_j = 2(N-1) + 2d$ , then $X \preccurlyeq \Delta_{S_d}$ implies that X is c-graphical.

Proof. This follows from Theorem 1 and the fact that all $S_d$ are connected graphs.

Comments

The case d = 0 is Hakimi's theorem. The main advantage of this result is that it is easier to verify if $X \preccurlyeq \Delta_{S_d}$ than to show directly that X is c-graphical. As the first component of $S_d$ is (N-1) it is relatively easy for a sequence X to be smaller than $\Delta_{S_d}$.

We note that $X \preccurlyeq \Delta_{S_d}$ is not a necessary condition for X to be c-graphical. Take N = 5 and X =(4,3,3,3,1). Then d = 3 and we know (Fig. 1) that X is c-graphical. Now, $\Delta_{S_3} = (4,4,2,2,2)$ and $X \not\preccurlyeq \Delta_{S_3}$ (and we note that also $\Delta_{S_3} \not\preccurlyeq X$).

We already know that X $\preccurlyeq \Delta_{S_d}$ implies that X is c-graphical and that the converse does not hold. Yet, we have the following necessary condition to be not graphical.

Notation. Given Y $\in \mathbb{N}^N = (y_1, \dots, y_N)$, then we denote by $Y_H$ the N-sequence $Y_H = \left(y_2 - 1, \dots, y_{y_1+1} - 1, y_{y_1+2}, \dots, y_N\right)$. The transformation that maps Y to $Y_H$ is called an H-transformation.

As a preparation of the next theorem we first note the following (trivial) lemma:

If X = $(x_1, \dots, x_N)$, and X'= $(x'_1, \dots, x'_N)$, then $x_1 = x'_1 = N - 1$, $\sum_{j=1}^{N} x_j = \sum_{j=1}^{N} x'_j$ , X $\prec X'$ implies $X_H \prec X'_H$. Indeed, in this case $X_H$=($x_2$-1,...,$x_N$-1) and $X'_H$=($x'_2$-1,...,$x'_N$-1) .□

This lemma is not true if $x_1$=$x'_1$, and they are not equal to N-1. Indeed, take N=5, X=(3,3,2,2,2) $\prec$ X'=(3,3,3,2,1). Then $X_H$=$X'_H$=(2,2,1,1).

Theorem 2

Given N, and $X \in \mathbb{N}^N$. If now $\sum_{j=1}^{N} x_j = 2(N-1) + 2d$ and $\Delta_{S_d} \prec X$, then X is not graphical.

Proof of the theorem. Consider $\Delta_{S_d}$ as described in equation (17) and apply the H-transformation i(d) times. The result is the following sequence of length N-i(d): $\left(j(d), \underbrace{1,\ldots,1}_{j(d) times}, \underbrace{0,\ldots 0}_{(N-i(d)-j(d)-1)\ times}\right)$. As $\Delta_{S_d} \prec X$ and the first i(d) coordinates of $\Delta_{S_d}$ are (N-1), this implies that also the first i(d) coordinates of X are (N-1). Applying now the H-transformation i(d) times (for N, N-1, …, N-i(d)+1, as sequence lengths) yields (by the lemma above):

$$\left(i(d), \underbrace{1,\ldots,1}_{j(d) times}, \underbrace{0,\ldots 0}_{(N-i(d)-j(d)-1)\ times}\right) \prec X_{\underbrace{H\ldots H}_{i(d) times}}$$

Now, $\left(i(d), \underbrace{1,\ldots,1}_{j(d) times}, \underbrace{0,\ldots 0}_{(N-i(d)-j(d)-1)\ times}\right)$ is the sequence corresponding to a star with i(d) nodes and (N-i(d)-j(d)-1) isolated points. This shows that $X_{\underbrace{H\ldots H}_{i(d) times}}$ is not graphical as a star has the highest delta sequence in the $\preccurlyeq$ sense. By the Havel-Hakimi theorem, it follows that X is not graphical.

An application: (4,4,3,2,1) is not graphical. Indeed, with N =5, $\Delta_{S_3}$= (4,4,2,2,2) $\prec$ (4,4,3,2,1), showing that (4,4,3,2,1) is not graphical.

We end this investigation by formulating a theoretical characterization for a sequence to be c-graphical.

Notation

The set of all maximal elements of the poset $\Delta(T^d(N))$, $\preccurlyeq$ is denoted as $M(T^d(N))$

Proposition 3

$X \in \mathbb{N}^N$ with $\sum_{j=1}^{N} x_j = 2(N-1) + 2d$ is c-graphical if and only if there exists $Y \in M(T^d(N))$ such that $X \preccurlyeq Y$.

Proof. One implication is an immediate consequence of Theorem 1, while the other implication follows from the fact that every element in a finite poset is smaller than or equal to a maximal element.

If we remove the requirement in $T^d(N)$ that the network is connected then we obtain a characterization for graphical networks, i.e., not necessarily connected ones.

If d=0, then $M(T^o(N)) = \{\left(N-1, \underbrace{1,\dots,1}_{(N-1)times}\right)\}$. Hence in this case we obtain Hakimi's theorem again.

This leads to the following open problem: determine $M(T^d(N))$ for all d ≥ 1. We will further provide a partial answer to this open problem.

Some further investigations related to $M(T^d(N))$

Theorem 3.

For all d ≥ 0: If $X = (x_1, \dots, x_N) \in M\left(T^d(N)\right)$, then $x_1 = N-1$

Proof. Assume that $x_1 < N-1$ then there are $N-1-x_1$ points which are not directly linked to the point 1. These points form a set, denoted as S. For every point j in A, there exists a path between j and 1 (as the network is assumed to be connected). We denote this path as j – $j_1$ - … - $j_n$ -1. We now remove the link between j and $j_1$ and add a link between j and 1. The new network is still connected. As X is decreasing we have $x_1 \geq x_{j_1} \geq 2$ (in the original network) and in the new network, with delta sequence denoted as $X' = (x'_1, \dots, x'_N)$, we have $x'_1 = x_1 + 1$, $x'_{j_1} = x_j - 1 \geq 1$, such that $X \prec X'$ (not equal!). This contradicts the fact that $X \in M\left(T^d(N)\right)$. Hence $x_1 = N-1$. □

Remark. This theorem does not hold for $x_2$ instead of $x_1$. We provide an example. Let N = 5, d=3 and consider X = (4,3,3,3,1). X is c-graphical as it is the degree sequence of Fig. 1.

Now, $x_2$ = 3 = N=2 < N-1 and the only sequence X' with $x'_2$ = N-1 and X≼X' (X≠X') is X'= (4,4,3,2,1), which is not graphical (and can easily be shown). Consequently $X \in M(T^3(N))$ and $x_2 \neq N-1$ (and hence $x_j \neq N-1, \text{for all}\, j \geq 2$).

We already know that for all d, $\Delta(S_d) \in M(T^d(N))$. Besides $S_d$ we further define $S'_d$ for all $d \in \left\{0,1.,,,,\frac{(N-1)(N-2)}{2}\right\}$. For d = 0,1,2 $S'_d = S_d$. For d ≥ 3, $S'_d$ is obtained, starting from the star $S_0(N)$ by constructing in the already obtained network a complete n-graph (n = 3,4,….). This is illustrated in Figs. 4 and 5.

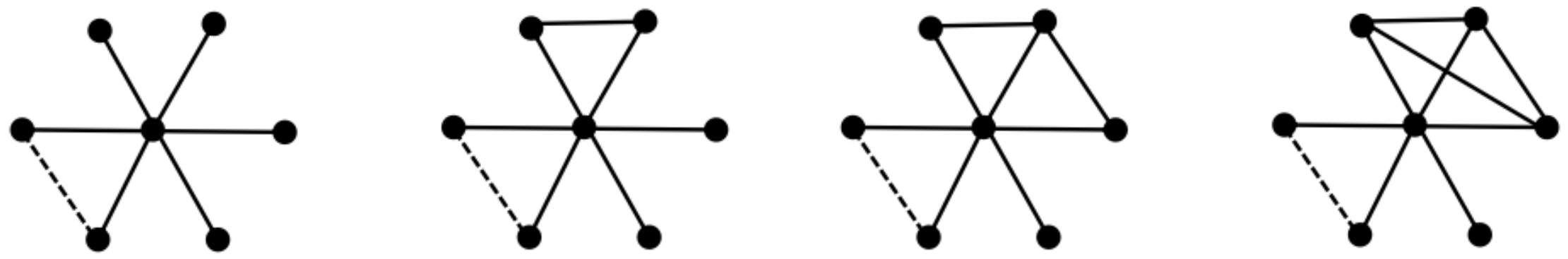

Fig.4. From left to right we have $S_0 = S'_0$, $S_1 = S'_1$, $S_2 = S'_2$, and $S'_3$ with $\Delta(S'_3) = \left(N-1,3,3,3,\ \underbrace{1,\dots,1}_{(N-4)times}\right)$

Next we have the cycle from n =3 to n=4, see Fig.5.

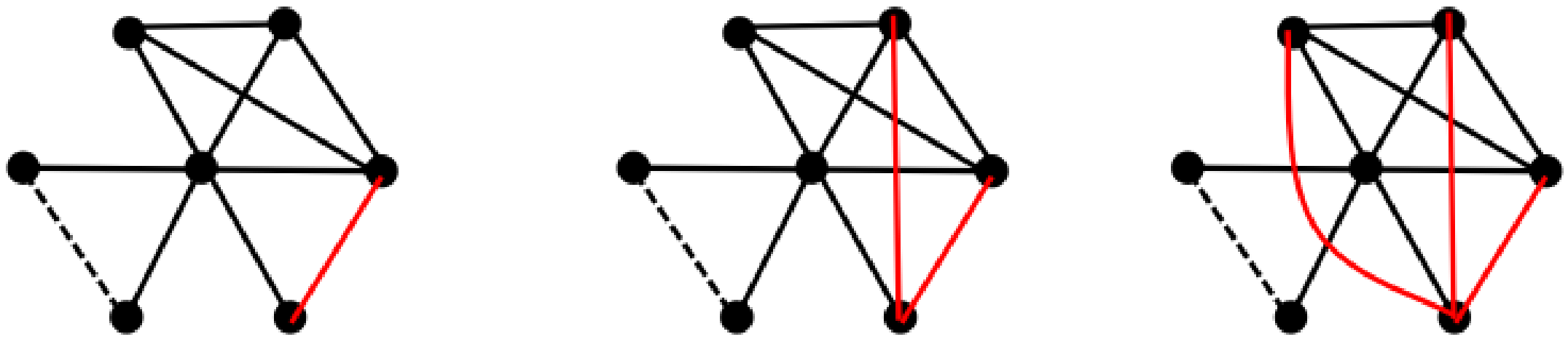

Fig. 5. From left to right we have: $S'_4, S'_5$ and $S'_6$.

We note that $\Delta(S'_4) = \left(N-1,4,3,3,2,\ \underbrace{1,\dots 1}_{(N-5)\ times}\right)$, $\Delta(S'_5) = \left(N-1,4,4,3,3,\ \underbrace{1,\dots 1}_{(N-5)\ times}\right)$, $\Delta(S'_6) = \left(N-1,4,4,4,4,\ \underbrace{1,\dots 1}_{(N-5)\ times}\right)$.

We will next determine $M(T^d(N))$ for d =0,1,2,3 and 4, which is a generalization of Hakimi's result (the case d=0). For this we need the following lemma.

Lemma

$$\forall d \geq 0: \Delta\left(S'_d\ (N)\right) \in M(T^D)$$

Proof. A general network $S'_d(N)$ consists of a completely connected n-graph (n = m(d)) and m diagonals (m = n(d)) of an (n+1) graph and points with δ =1. For to step to d+1 there are two possible cases. First, the (n+1) graph is completely connected. The only way to reach $S'_{d+1}(N)$ is by connecting a point with δ =1 to a vertex of this complete (n+1) graph. Second, the (n+1) graph is not yet completely connected. Connecting further to reach $S'_{d+1}(N)$ yields a Δ-sequence which is higher in the ≼ partial order than the sequence obtained by connecting a point with δ = 1 with the (n=1) graph which was not yet completely connected. This proves this lemma. □

An illustration is given in Fig. 6.

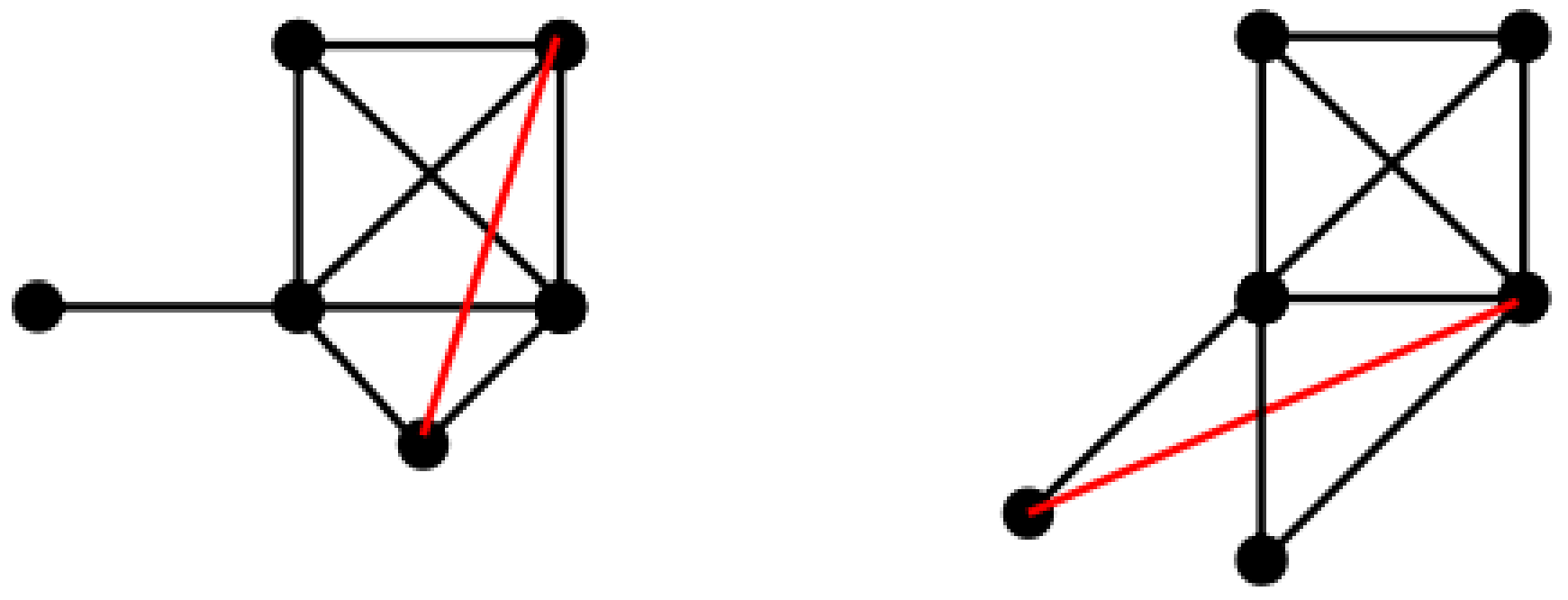

Fig. 6. The black part is $S'_d$. The left-hand figure is $S'_{d+1}$ and the right-hand one is S such that $\Delta(S) \preccurlyeq \Delta(S'_{d+1})$.

We note that we already know that $\Delta(S_d(N)) \in M(T^d)$, hence by this lemma we have $\{\Delta(S_d(N)), \Delta(S'_d(N))\} \subset M(T^d)$.

Now we come to the characterization theorem.

Theorem 4. Extension of Theorem of Hakimi

$\forall X = (x_1, \dots, x_N) \in \mathbb{N}^N$, such that $\sum_{j=1}^{N} x_j$ = 2(N-1)+d, we have

1) If d=0 (Hakimi): M($T^0$) = $\{\Delta(S_0(N))\} = \left\{(N-1, \underbrace{1,\dots,1}_{(N-1)times})\right\}$. X is c-graphical, representing a tree.
2) If d= 1:M($T^1$) = $\{\Delta(S_1(N))\} = \left\{(N-1,2,2, \underbrace{1,\dots,1}_{(N-3)times})\right\}$ and X is c-graphical if and only if X $\preccurlyeq$ $\Delta(S_1(N))$ .
3) If d =2: M($T^2$) = $\{\Delta(S_2(N))\} = \left\{(N-1,3,2,2, \underbrace{1,\dots,1}_{(N-4)times})\right\}$ and X is c-graphical if and only if X $\preccurlyeq$ $\Delta(S_2(N))$
4) If d = 3: : M($T^3$) = $\{\Delta(S_3(N)), \Delta(S_3'(N))\}$ = $\left\{\left(N-1,4,2,2,2, \underbrace{1,\dots,1}_{(N-5)times}\right), \left(N-1,3,3,3,, \underbrace{1,\dots,1}_{(N-4)times}\right)\right\}$ and X is c-graphical if and only if (X $\preccurlyeq$ $\Delta(S_3(N))$ OR X $\preccurlyeq$ $\Delta(S_3'(N))$
5) If d = 4: M($T^4$) = $\{\Delta(S_4(N)), \Delta(S_4'(N))\}$ = $\left\{\left(N-1,5,2,2,2,2, \underbrace{1,\dots,1}_{(N-6)times}\right), \left(N-1,4,3,3,2, \underbrace{1,\dots,1}_{(N-5)times}\right)\right\}$ and X is c-graphical if and only if (X $\preccurlyeq$ $\Delta(S_4(N))$ OR X $\preccurlyeq$ $\Delta(S_4'(N))$
6) If d ≥ 5: M($T^d$) $\supsetneq$ $\{\Delta(S_d(N)),\ \Delta(S_d'(N))\}$.

Proof. By theorem 3 we know that we only have to consider sequences X such that $x_1$ = N-1.

1) d=0. As $\Delta(S_0(N)) = \left(N-1, \underbrace{1,\dots,1}_{(N-1)times}\right)$ is the largest in the poset of $\preccurlyeq$, every X such that $\sum_{j=1}^{N} x_j$ = 2(N-1) is c-graphical (Proposition 3). It is moreover well-known that a c-graphical X with $\sum_{j=1}^{N} x_j$ = 2(N-1) is a tree (West, 2001).
2) d=1. Now $S_1(N) = S_1'(N)$ and hence $\Delta(S_1(N)) = \Delta(S_1'(N))$ . Only $X = \left(N-1,3, \underbrace{1,\dots,1}_{(N-2)times}\right)$ satisfies the relation is strictly larger than $\Delta(S_1(N))$ in the Lorenz majorization order and hence, by Theorem 2, X is not graphical. Hence $M(T^1) = \{\Delta(S_1(N))\}$, from which the result follows by Proposition 3.
3) d=2. As for the case d=1, $S_2(N) = S_2'(N)$ and hence $\Delta(S_2(N)) = \Delta(S_2'(N))$. Now for d=2, there exists exactly one graphical variant S of $S_2$, namely the network shown in Fig. 7.

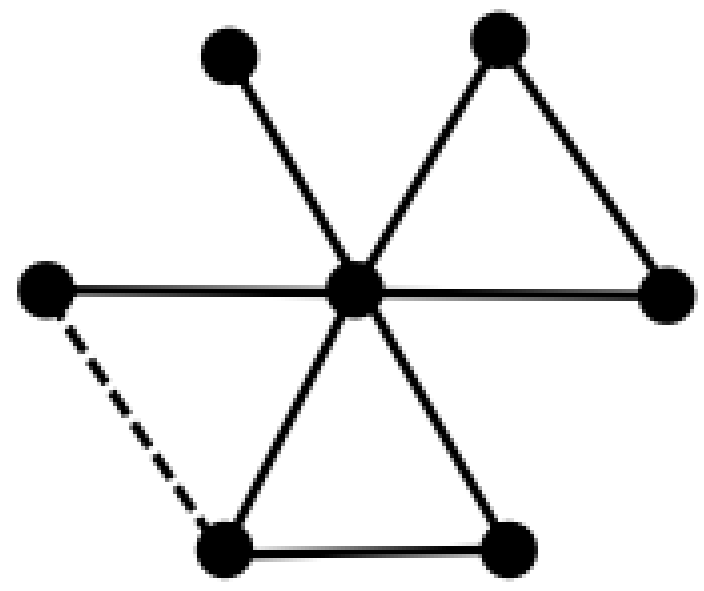

Fig. 7

with $\Delta(S) = \left(N-1,2,2,2,2, \underbrace{1,\dots,1}_{(N-5)times}\right) \preccurlyeq \Delta(S_2)$ (strictly), by Theorem 3. Hence $\Delta(S) \notin M(T^2)$ and thus *M(T²)* = $\{\Delta(S_2(N))\}$, leading to the required result by Proposition 3.
4) d=3. By the previous proof, we only have to consider *$S_2(N)$* and find out how we can add – maximally – one more link

(d=3). This leads to $S_2(N), S_2'(N)$ or the graph S shown in Fig. 8.

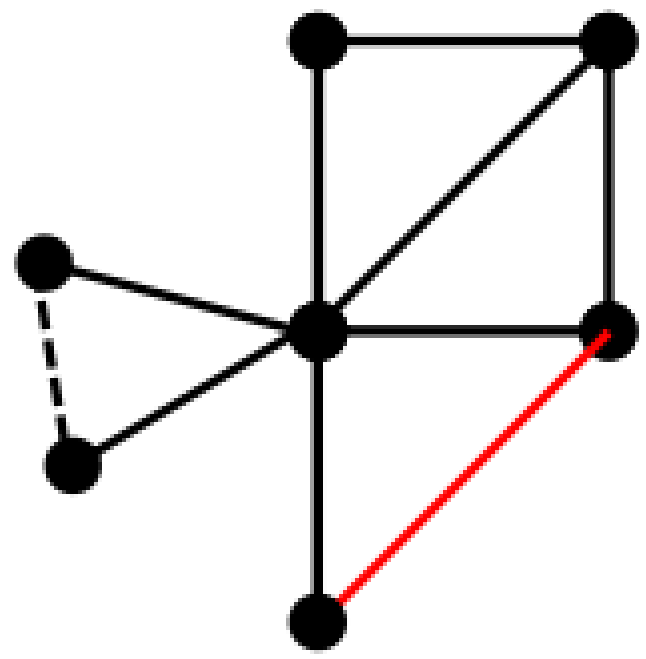

Fig. 8

with $\Delta(S(N)) = (N-1, 3{,}3{,}2{,}2 \underbrace{,1, \dots ,1)}_{(N-5) times} \preccurlyeq \Delta(S_3'(N))$ . Hence $\Delta(S(N)) \notin M(T^3)$ . As, moreover, $\Delta(S_3(N))$ and $\Delta(S_3'(N))$ are incomparable in the Lorenz majorization poset we have that $M(T^3) = \{\Delta(S_3(N)), \Delta(S_3'(N))\}$ and the result follows from Proposition 3.

5) d=4. By part 4) we need only investigate $S_3$ and $S_3'$ and see how we can add one more link (in a maximal fashion) (d=4). Via $S_3'$ this can only as $S_4'$ and S (see Fig. 9) with $\Delta(S(N)) = \left(N-1,3,3,3,2,2,\underbrace{1,\dots,1}_{(N-6)times}\right) \preccurlyeq \Delta\left(S_3'(N)\right)$. Hence $\Delta(S(N)) \notin M(T^4)$, or, via $S_3$, as $S_4$ and one of the three graphs $S^{(k)}, k = 1,2,3$, shown in Fig. 10.

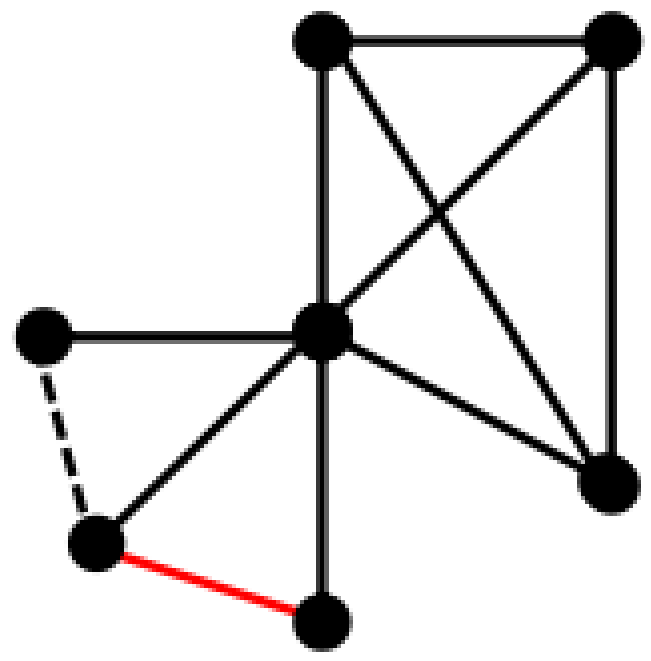

Fig. 9

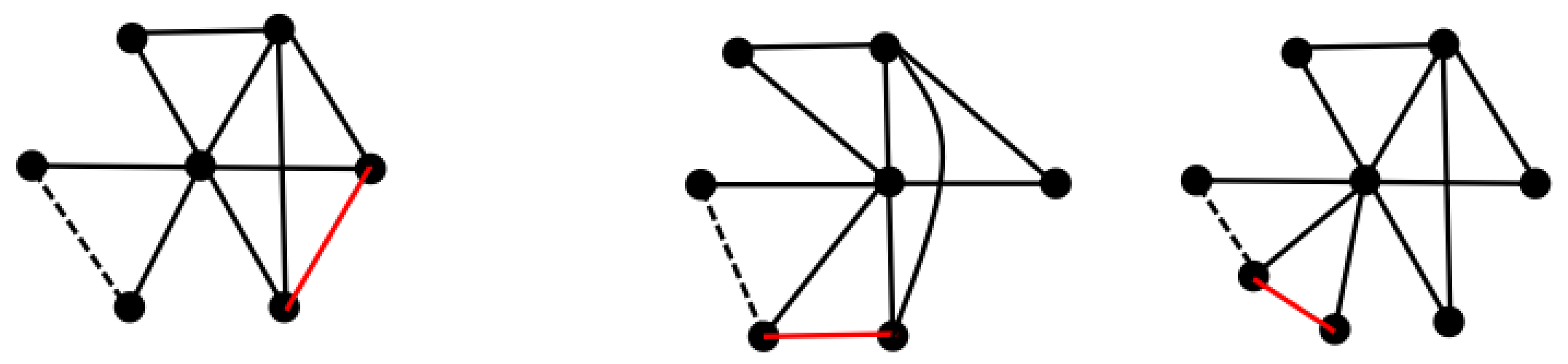

Fig. 10. This figure shows the three graphs $S^{(k)}, k = 1,2,3$

We see that $\Delta\left(S^{(1)}(N)\right) = \left(N-1,4,3,3,2,\underbrace{1,\dots,1}_{(N-5)times}\right) = \Delta(S_4'(N))$;

$\Delta\left(S^{(2)}(N)\right) = \left(N-1,4,3,2,2,2,\underbrace{1,\dots,1}_{(N-6)times}\right) \preccurlyeq \Delta(S_4'(N))$ ;

$\Delta\left(S^{(3)}(N)\right) = \left(N-1,4,2,2,2,2,2,\underbrace{1,\dots,1}_{(N-7)times}\right) \preccurlyeq \Delta(S_4'(N))$ . This shows that in each case $\Delta\left(S^{(k)}(N)\right) \notin M(T^4)$. As moreover, $\Delta(S_4(N))$ and $\Delta(S_4'(N))$ are incomparable in the Lorenz majorization poset we have that $M(T^4) = \{\Delta(S_4(N)), \Delta(S_4'(N))\}$ and the result follows from Proposition 3.

6) Finally, we only have to find a network *S(N)* which is incomparable with $\Delta(S_5(N))$, as well as with $\Delta(S_5'(N))$. Such a network S is shown in Fig. 11.

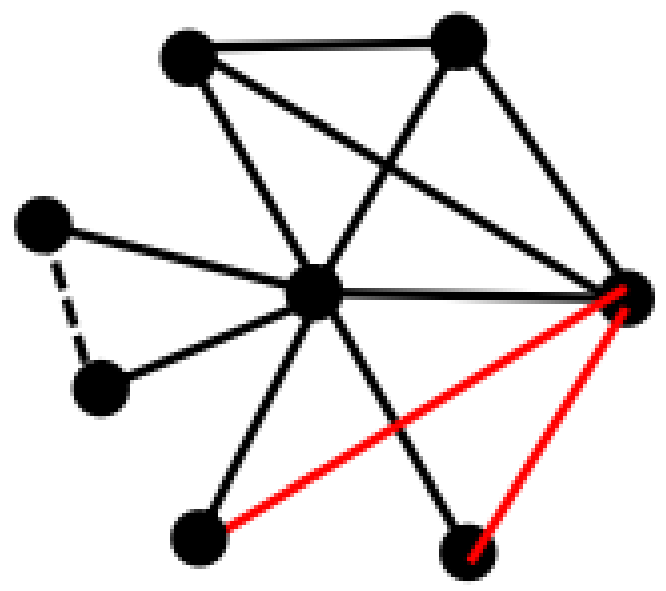

Fig. 11. The network *S* needed in point 6).

We have: $\Delta(S) = \left(N-1, 5,3,3,2,2, \underbrace{1,\dots,1}_{(N-6)times}\right)$. It is now easy to see that $\Delta(S)$ is incomparable with $\Delta(S_5(N)) = \left(N-1, 6,2,2,2,2,2, \underbrace{1,\dots,1}_{(N-7)times}\right)$ as well as with $\Delta(S_5'(N)) = \left(N-1, 4,4,3,3 \underbrace{1,\dots,1}_{(N-5)times}\right)$. This shows that $M(T^5) \supsetneq \{\Delta(S_5(N)),\ \Delta(S_5'(N))\}$. This same idea can be used to prove that, for all d ≥5: $M(T^d) \supsetneq \{\Delta(S_d(N)),\ \Delta(S_d'(N))\}$.□

The case d < 0

All networks with degree sequence $X \in \mathbb{N}^N$ satisfy the relation that the sum of all degrees is larger than or equal to zero, where the sum is equal to zero for a set of N unconnected points. Considering the formula $\sum_{j=1}^{N} x_j = 2(N-1) + 2d$ shows that the case $d \in \{-(N-1), -(N-2), \dots, -1\}$ also exists. In this case $S_d(N)$ denotes the incomplete star constructed by connecting one point with $\frac{\sum_{j=1}^{N} x_j}{2} = N - 1 + d < (N-1)$ points. The case d < 0 can – trivially – be added to Hakimi's existence theorem (Theorem 4):

Theorem 5. For all sequences X as described above (d < 0) X is graphical.

Proof. We see that $\Delta(S_d(N)) = \left(a, \underbrace{1, \ldots, 1}_{a\ times}, \underbrace{0, \ldots, 0}_{(N-a-1)\ times}\right)$ with 0 ≤ a < N-1, and hence for each network with $\sum_{j=1}^{N} x_j = 2a$,we have $X \preccurlyeq S_d(N)$. By Theorem 1 it follows that X is graphical (not necessarily c-graphical as $S_d(N)$ is not connected).□

## Conclusions

This paper studies how the Lorenz majorization order on degree sequences of a network can determine whether or not a degree sequence is graphical (or not), or even c-graphical (or not). The basic idea behind this is that it is much easier to determine Lorenz majorization than whether or not a finite sequence is graphical. The main result states that if X is dominated by X' (in the Lorenz order) and if X' is graphical, respectively c-graphical, then X is also graphical (respectively c-graphical). The proof is direct in the sense that it does not use existing realizability theorems, such as Erdös-Gallai or Havel-Hakimi. Yet we present a second shorter proof of our result, which uses the former theorem.

From this general result, we reprove and extend the famous Havel-Hakimi theorem. Moreover, not only does the classical theorem of Hakimi about trees follow from this general result, but also a generalization to general connected networks. We further prove a sufficient condition for a degree sequence not to be graphical (this result uses Havel-Hakimi) and leads to a characterization of c-graphical sequences: a degree sequence X of length N with sum equal to 2(N-1) +2d is c-graphical if and only if X is dominated (for the Lorenz order) by a maximal finite sequence of equal length, in the space of all degree sequences of connected graphs with sum equal to 2(N-1) + 2d.

It would be interesting to give concrete expressions of the maximal finite sequences. The paper closes with some examples for this, for every d.

Acknowledgements. The author thanks Ronald Rousseau for useful discussions and Li Li for providing the figures.

No funding was obtained for this investigation.

This paper is a revision and extension of (Egghe, 2024).